\documentclass[12pt]{article}
\usepackage{amsfonts}

\usepackage{myart}

\normalbaselines
\input tcilatex

\begin{document}

\title{Coordinateless description and deformation principle as foundations of
physical geometry}
\author{Yuri A.Rylov}
\date{Institute for Problems in Mechanics, Russian Academy of Sciences,\\
101-1, Vernadskii Ave., Moscow, 119526, Russia.\\
e-mail: rylov@ipmnet.ru\\
Web site: {$http://rsfq1.physics.sunysb.edu/\symbol{126}rylov/yrylov.htm$}\\
or mirror Web site: {$http://gasdyn-ipm.ipmnet.ru/\symbol{126}%
rylov/yrylov.htm$}}
\maketitle

\begin{abstract}
Physical geometry studies mutual disposition of geometrical objects and
points in space, or space-time, which is described by the distance function $%
d$, or by the world function $\sigma =d^{2}/2$. One suggests a new general
method of the physical geometry construction. The proper Euclidean geometry
is described in terms of its world function $\sigma _{\mathrm{E}}$. Any
physical geometry $\mathcal{G}$ is obtained from the Euclidean geometry as a
result of replacement of the Euclidean world function $\sigma _{\mathrm{E}}$
by the world function $\sigma $ of $\mathcal{G}$. This method is very simple
and effective. It introduces a new geometric property: nondegeneracy of
geometry. Using this method, one can construct deterministic space-time
geometries with primordially stochastic motion of free particles and
geometrized particle mass. Such a space-time geometry defined properly (with
quantum constant as an attribute of geometry) allows one to explain quantum
effects as a result of the statistical description of the stochastic
particle motion (without a use of quantum principles).
\end{abstract}

\section{Introduction}

A geometry lies in the foundation of physics, and a true conception of
geometry is very important for the consequent development of physics. It is
common practice to think that all problems in foundations of geometry have
been solved many years ago. It is valid, but this concerns the geometry
considered to be a logical construction. Physicists are interested in the
geometry considered as a science on mutual disposition of geometrical
objects in the space or in the space-time. The two aspects of geometry are
quite different, and one can speak about two different geometries, using for
them two different terms. Geometry as a logical construction is a
homogeneous geometry, where all points have the same properties. Well known
mathematician Felix Klein \cite{K37} believed that only the homogeneous
geometry deserves to be called a geometry. It is his opinion that the
Riemannian geometry (in general, inhomogeneous geometry) should be qualified
as a Riemannian geography, or a Riemannian topography. In other words, Felix
Klein considered a geometry mainly as a logical construction. We shall refer
to such a geometry as the mathematical geometry.

The geometry considered to be a science on mutual disposition of geometric
objects will be referred to as a physical geometry, because the physicists
are interested mainly in this aspect of a geometry. The physical geometries
are inhomogeneous, in general, although they may be homogeneous also. On the
one hand, the proper Euclidean geometry is a physical geometry. On the other
hand, it is a logical construction, because it is homogeneous and can be
constructed of simple elements (points, straights, planes, etc.). All
elements of the Euclidean geometry have similar properties, which are
described by axioms. Similarity of geometrical elements allows one to
construct the mathematical (homogeneous) geometry by means of logical
reasonings. The proper Euclidean geometry was constructed many years ago by
Euclid. Consistency of this construction was investigated and proved in \cite
{H30}. Such a construction is very complicated even in the case of the
proper Euclidean geometry, because simple geometrical objects are used for
construction of more complicated ones, and one cannot construct a
complicated geometrical object $\mathcal{O}$ without construction of more
simple constituents of this object.

Note that constructing his geometry, Euclid did not use coordinates for
labeling of the space points. His description of the homogeneous geometry
was coordinateless. It means that the coordinates are not a necessary
attribute of the geometry. Coordinate system is a method of the geometry
description, which may or may not be used. Application of coordinates and of
other means of description poses the problem of separation of the geometry
properties from the properties of the means of the description. Usually the
separation of the geometry properties from the coordinate system properties
is carried out as follows. The geometry is described in all possible
coordinate systems. Transformations from one coordinate system to the
another one form a group of transformation. Invariants of this
transformation group are the same in all coordinate system, and hence, they
describe properties of the geometry in question.

At this point we are to make a very important remark. \textit{Any geometry
is a totality of all geometric objects }$\mathcal{O}$\textit{\ and of all
relations }$\mathcal{R}$\textit{\ between them}. Any geometric object $%
\mathcal{O}$ is a subset of points of some point set $\Omega $, where the
geometry is given. In the Riemannian geometry (and in other inhomogeneous
geometries) the set $\Omega $ is supposed to be a $n$-dimensional manifold,
whose points $P$ are labelled by $n$ coordinates $x=\left\{
x^{1},x^{2},...x^{n}\right\} $. This labelling (arithmetization of space) is
considered to be a necessary attribute of the Riemannian geometry. Most
geometers believe that the Riemannian geometry (and physical geometry, in
general) cannot be constructed without introduction of the manifold. In
other words, they believe that the manifold is an attribute of the
Riemannian geometry (and of any continuous geometry, in general). This
belief is founded on the fact, that the Riemannian geometry is always
constructed on some manifold. But this belief is a delusion. The fact, that
we always construct the physical geometry on some manifold, does not mean
that the physical geometry cannot be constructed without a reference to a
manifold, or to a coordinate system. Of course, some labelling of the
spatial points (coordinate system) is convenient, but this labelling has no
relations to the construction of the geometry, and the physical geometry
should be constructed without a reference to coordinate system. Application
of the coordinate system imposes constraints on properties of the
constructed physical geometry. For instance, if we use a continuous
coordinate system (manifold) we can construct only continuous physical
geometry. To construct a discrete physical geometry, the geometry
construction is not to contain a reference to the coordinate system.

Here we present the method of the physical geometry construction, which does
not contain a reference to the coordinate system and other means of
description. It contains a reference only to the distance function $d$,
which is a real characteristic of physical geometry.

If a geometry is inhomogeneous, and the straights located in different
places have different properties, it is impossible to describe properties of
straights by means of axioms, because there are no such axioms for the whole
geometry. Mutual disposition of points in a physical (inhomogeneous)
geometry, which is given on the set $\Omega $ of points $P$, is described by
the distance function $d\left( P,Q\right) $%
\begin{equation}
d:\qquad \Omega \times \Omega \rightarrow \Bbb{R},\qquad d\left( P,P\right)
=0,\qquad \forall P\in \Omega  \label{b1.1}
\end{equation}
The distance function $d$ is the main characteristic of the physical
geometry. Besides, the distance function $d$ is an \textit{unique
characteristic} of any physical geometry. \textit{The distance function }$d$%
\textit{\ determines completely the physical geometry.} This statement is
very important for construction of a physical geometry. It will be proved
below. Any physical geometry $\mathcal{G}$ is constructed from the proper
Euclidean geometry $\mathcal{G}_{\mathrm{E}}$ by means of a deformation,
i.e. by a replacement of the Euclidean distance function $d_{\mathrm{E}}$ by
the distance function of the geometry in question. For instance,
constructing the Riemannian geometry, we replace the Euclidean infinitesimal
distance $dS_{\mathrm{E}}=\sqrt{g_{\mathrm{E}ik}dx^{i}dx^{k}}$ by the
Riemannian one $dS=\sqrt{g_{ik}dx^{i}dx^{k}}$. There is no method of the
inhomogeneous physical geometry construction other, than the deformation of
the Euclidean geometry (or some other homogeneous geometry) which is
constructed as a mathematical geometry on the basis of its axiomatics and
logic. Unfortunately, conventional method of the Riemannian geometry
construction contains a reference to the coordinate system. But this
reference can be eliminated, provided that we use finite distances $d$
instead of infinitesimal distances $dS$.

For description of a physical geometry one uses the world function $\sigma $ 
\cite{S60}, which is connected with the distance function $d$ by means of
the relation $\sigma \left( P,Q\right) =\frac{1}{2}d^{2}\left( P,Q\right) $.
The world function $\sigma $ of the $\sigma $-space $V=\left\{ \sigma
,\Omega \right\} $ is defined by the relation 
\begin{equation}
\sigma :\qquad \Omega \times \Omega \rightarrow \Bbb{R},\qquad \sigma \left(
P,P\right) =0,\qquad \forall P\in \Omega  \label{b1.2}
\end{equation}
Application of the world function is more convenient in the relation that
the world function is real, when the distance function $d$ is imaginary and
does not satisfy definition (\ref{b1.1}). It is important at the
consideration of the space-time geometry as a physical geometry.

In general, a physical geometry cannot be constructed as a logical building,
because any change of the world function should be accompanied by a change
of axiomatics. This is practically aerial, because the set of possible
physical geometries is continual. Does the world function contain full
information which is necessary for construction of the physical geometry? It
is a very important question. For instance, can one derive the space
dimension from the world function in the case of Euclidean geometry?
Slightly below we shall answer this question in the affirmative. Now we
formulate the method of the physical geometry construction.

Let us imagine that the proper Euclidean geometry $\mathcal{G}_{\mathrm{E}}$
can be described completely in terms and only in terms of the Euclidean
world function $\sigma _{\mathrm{E}}$. Such a description is called $\sigma $%
-immanent. It means that any geometrical object $\mathcal{O}_{\mathrm{E}}$
and any relation $\mathcal{R}_{\mathrm{E}}$ between geometrical objects in $%
\mathcal{G}_{\mathrm{E}}$ can be described in terms of $\sigma _{\mathrm{E}}$
in the form $\mathcal{O}_{\mathrm{E}}\left( \sigma _{\mathrm{E}}\right) $
and $\mathcal{R}_{\mathrm{E}}\left( \sigma _{\mathrm{E}}\right) $. To obtain
corresponding geometrical object $\mathcal{O}$ and corresponding relation $%
\mathcal{R}$ between the geometrical objects in other physical geometry $%
\mathcal{G}$, it is sufficient to replace the Euclidean world function $%
\sigma _{\mathrm{E}}$ by the world function $\sigma $ of the physical
geometry $\mathcal{G}$ in description of $\mathcal{O}_{\mathrm{E}}\left(
\sigma _{\mathrm{E}}\right) $ and $\mathcal{R}_{\mathrm{E}}\left( \sigma _{%
\mathrm{E}}\right) $. 
\[
\mathcal{O}_{\mathrm{E}}\left( \sigma _{\mathrm{E}}\right) \rightarrow 
\mathcal{O}_{\mathrm{E}}\left( \sigma \right) ,\qquad \mathcal{R}_{\mathrm{E}%
}\left( \sigma _{\mathrm{E}}\right) \rightarrow \mathcal{R}_{\mathrm{E}%
}\left( \sigma \right) 
\]
Index 'E' shows that the geometric object is constructed on the basis of the
Euclidean axiomatics. Thus, one can obtain another physical geometry $%
\mathcal{G}$ from the Euclidean geometry $\mathcal{G}_{\mathrm{E}}$ by a
simple replacement of $\sigma _{\mathrm{E}}$ by $\sigma $. For such a
construction one needs no axiomatics and no reasonings. One needs no means
of descriptions (topological structures, continuity, coordinate system,
manifold, dimension, etc.). In fact, one uses implicitly the axiomatics of
the Euclidean geometry, which is deformed by the replacement $\sigma _{%
\mathrm{E}}\rightarrow \sigma $. This replacement may be interpreted as a
deformation of the Euclidean space. Absence of a reference to the means of
description is an advantage of the considered method of the geometry
construction. Besides, there is no necessity to construct the whole geometry 
$\mathcal{G}$. We can construct and investigate only that part of the
geometry $\mathcal{G}$ which we are interested in. Any physical geometry may
be constructed as a result of a deformation of the Euclidean geometry.

The geometric object $\mathcal{O}$ is described by means of the
skeleton-envelope method \cite{R01}. It means that any geometric object $%
\mathcal{O}$ is considered to be a set of intersections and joins of
elementary geometric objects (EGO).

The finite set $\mathcal{P}^{n}\equiv \left\{ P_{0},P_{1},...,P_{n}\right\}
\subset \Omega $ of parameters of the envelope function $f_{\mathcal{P}^{n}}$
is the skeleton of elementary geometric object (EGO) $\mathcal{E}\subset
\Omega $. The set $\mathcal{E}\subset \Omega $ of points forming EGO is
called the envelope of its skeleton $\mathcal{P}^{n}$. For continuous
physical geometry the envelope $\mathcal{E}$ is usually a continual set of
points. The envelope function $f_{\mathcal{P}^{n}}$, determining EGO is a
function of the running point $R\in \Omega $ and of parameters $\mathcal{P}%
^{n}\subset \Omega $. The envelope function $f_{\mathcal{P}^{n}}$ is
supposed to be an algebraic function of $s$ arguments $w=\left\{
w_{1},w_{2},...w_{s}\right\} $, $s=(n+2)(n+1)/2$. Each of arguments $%
w_{k}=\sigma \left( Q_{k},L_{k}\right) $ is a $\sigma $-function of two
arguments $Q_{k},L_{k}\in \left\{ R,\mathcal{P}^{n}\right\} $, either
belonging to skeleton $\mathcal{P}^{n}$, or coinciding with the running
point $R$. Thus, any elementary geometric object $\mathcal{E}$ is determined
by its skeleton and its envelope function.

For instance, the sphere $\mathcal{S}(P_{0},P_{1})$ with the center at the
point $P_{0}$ is determined by the relation 
\begin{equation}
\mathcal{S}(P_{0},P_{1})=\left\{ R|f_{P_{0}P_{1}}\left( R\right) =0\right\}
,\qquad f_{P_{0}P_{1}}\left( R\right) =\sqrt{2\sigma \left(
P_{0},P_{1}\right) }-\sqrt{2\sigma \left( P_{0},R\right) }  \label{b1.6}
\end{equation}
where $P_{1}$ is a point belonging to the sphere. The elementary object $%
\mathcal{E}$ is determined in all physical geometries at once. In
particular, it is determined in the proper Euclidean geometry, where we can
obtain its meaning. We interpret the elementary geometrical object $\mathcal{%
E}$, using our knowledge of the proper Euclidean geometry. Thus, the proper
Euclidean geometry is used as a sample geometry for interpretation of any
physical geometry.

We do not try to repeat subscriptions of Euclid at construction of the
geometry. We take the geometrical objects and relations between them,
prepared in the framework of the Euclidean geometry and describe them in
terms of the world function. Thereafter we deform them, replacing the
Euclidean world function $\sigma _{\mathrm{E}}$ by the world function $%
\sigma $ of the geometry in question. In practice the construction of the
elementary geometry object is reduced to the representation of the
corresponding Euclidean geometrical object in the $\sigma $-immanent form,
i.e. in terms of the Euclidean world function. The last problem is the
problem of the proper Euclidean geometry. The problem of representation of
the geometrical object (or relation between objects) in the $\sigma $%
-immanent form is a real problem of the physical geometry construction.

It is very important, that such a construction does not use coordinates and
other methods of description, because the application of the means of
description imposes constraints on the constructed geometry. Any means of
description is a structure $St$ given on the basic Euclidean geometry with
the world function $\sigma _{\mathrm{E}}$. Replacement $\sigma _{\mathrm{E}%
}\rightarrow \sigma $ is sufficient for construction of unique physical
geometry $\mathcal{G}_{\sigma }$. If we use an additional structure $St$ for
construction of physical geometry, we obtain, in general, other geometry $%
\mathcal{G}_{St}$, which coincide with $\mathcal{G}_{\sigma }$ not for all $%
\sigma $, but only for some of world functions $\sigma $. Thus, a use of
additional means of description restricts the list of possible physical
geometries. For instance, if we use the coordinate description at
construction of the physical geometry, the obtained geometry appears to be
continuous, because description by means of the coordinates is effective
only for continuous geometries, where the number of coordinates coincides
with the geometry dimension.

Constructing geometry $\mathcal{G}$ by vtans of a deformation we use
essentially the fact that the proper Euclidean geometry $\mathcal{G}_{%
\mathrm{E}}$ is a mathematical geometry, which has been constucted on the
basis of Euclidean axiomatics and logical reasonings.

We shall refer to the described method of the physical geometry construction
as the deformation principle and interpret the deformation in the broad
sense of the word. In particular, a deformation of the Euclidean space may
transform an Euclidean surface into a point, and an Euclidean point into a
surface. Such a deformation may remove some points of the Euclidean space,
violating its continuity, or decreasing its dimension. Such a deformation
may add supplemental points to the Euclidean space, increasing its
dimension. In other words, the deformation principle is a very general
method of the physical geometry construction.

The deformation principle as a method of the physical geometry construction
contains two essential stages:

(i) Representation of geometrical objects $\mathcal{O}$ and relations $%
\mathcal{R}$ of the Euclidean geometry in the $\sigma $-immanent form, i.e.
in terms and only in terms of the world function $\sigma _{\mathrm{E}}$.

(ii) Replacement of the Euclidean world function $\sigma _{\mathrm{E}}$ by
the world function $\sigma $ of the geometry in question.

A physical geometry, constructed by means of the only deformation principle
(i.e. without a use of other methods of the geometry construction) is called
T-geometry (tubular geometry) \cite{R90,R01,R002}. The T-geometry is the
most general kind of the physical geometry.

Application of the deformation principle is restricted by two constraints.

1. Describing Euclidean geometric objects $\mathcal{O}\left( \sigma _{%
\mathrm{E}}\right) $ and Euclidean relation $\mathcal{R}\left( \sigma _{%
\mathrm{E}}\right) $ in terms of $\sigma _{\mathrm{E}}$, we are not to use
special properties of Euclidean world function $\sigma _{\mathrm{E}}$. In
particular, definitions of $\mathcal{O}\left( \sigma _{\mathrm{E}}\right) $
and $\mathcal{R}\left( \sigma _{\mathrm{E}}\right) $ are to have similar
form in Euclidean geometries of different dimensions. They must not depend
on the dimension of the Euclidean space.

2. The deformation principle is to be applied separately from other methods
of the geometry construction. In particular, one may not use topological
structures in construction of a physical geometry, because for effective
application of the deformation principle the obtained physical geometry must
be determined only by the world function (metric).

\section{Description of the proper Euclidean space in terms of the world
function}

The crucial point of the T-geometry construction is the description of the
proper Euclidean geometry in terms of the Euclidean world function $\sigma _{%
\mathrm{E}}$. We shall refer to this method of description as the $\sigma $%
-immanent description. Unfortunately, it was unknown for many years,
although all physicists knew that the infinitesimal interval $dS=\sqrt{%
g_{ik}dx^{i}dx^{k}}$ is the unique essential characteristic of the
space-time geometry, and changing this expression, we change the space-time
geometry. From physical viewpoint the $\sigma $-immanent description is very
reasonable, because it does not contain any extrinsic information. The $%
\sigma $-immanent description does not refer to the means of description
(dimension, manifold, coordinate system). Absence of references to means of
description is important in the relation, that there is no necessity to
separate the information on the geometry in itself from the information on
the means of description. The $\sigma $-immanent description contains only
essential characteristic of geometry: its world function. At first the $%
\sigma $-immanent description was obtained in 1990 \cite{R90}.

The first question concerning the $\sigma $-immanent description is as
follows. Does the world function contain sufficient information for
description of a physical geometry? The answer is affirmative, at least, in
the case of the proper Euclidean geometry, and this answer is given by the
prove of the following theorem.

Let $\sigma $-space $V=\left\{ \sigma ,\Omega \right\} $ be a set $\Omega $
of points $P$ with the given world function $\sigma $%
\begin{equation}
\sigma :\qquad \Omega \times \Omega \rightarrow \Bbb{R},\qquad \sigma \left(
P,P\right) =0,\qquad \forall P\in \Omega  \label{a1.3}
\end{equation}
Let the vector $\mathbf{P}_{0}\mathbf{P}_{1}\mathbf{=}\left\{
P_{0},P_{1}\right\} $ be the ordered set of two points $P_{0}$, $P_{1}$, and
its length $\left| \mathbf{P}_{0}\mathbf{P}_{1}\right| $ is defined by the
relation $\left| \mathbf{P}_{0}\mathbf{P}_{1}\right| ^{2}=2\sigma \left(
P_{0},P_{1}\right) $.

\textit{Theorem}

The $\sigma $-space $V=\left\{ \sigma ,\Omega \right\} $ is the $n$%
-dimensional proper Euclidean space, if and only if the world function $%
\sigma $ satisfies the following conditions, written in terms of the world
function $\sigma $.

I. Condition of symmetry: 
\begin{equation}
\sigma \left( P,Q\right) =\sigma \left( Q,P\right) ,\qquad \forall P,Q\in
\Omega  \label{a1.4}
\end{equation}

II. Definition of the dimension: 
\begin{equation}
\exists \mathcal{P}^{n}\equiv \left\{ P_{0},P_{1},...P_{n}\right\} ,\qquad
F_{n}\left( \mathcal{P}^{n}\right) \neq 0,\qquad F_{k}\left( {\Omega }%
^{k+1}\right) =0,\qquad k>n  \label{b10}
\end{equation}
where $F_{n}\left( \mathcal{P}^{n}\right) $ is the Gram's determinant 
\begin{equation}
F_{n}\left( \mathcal{P}^{n}\right) =\det \left| \left| \left( \mathbf{P}_{0}%
\mathbf{P}_{i}.\mathbf{P}_{0}\mathbf{P}_{k}\right) \right| \right| =\det
\left| \left| g_{ik}\left( \mathcal{P}^{n}\right) \right| \right| ,\qquad
i,k=1,2,...n  \label{b11}
\end{equation}
The scalar product $\left( \mathbf{P}_{0}\mathbf{P}_{1}.\mathbf{Q}_{0}%
\mathbf{Q}_{1}\right) $ of two vectors $\mathbf{P}_{0}\mathbf{P}_{1}$ and $%
\mathbf{Q}_{0}\mathbf{Q}_{1}$ is defined by the relation 
\begin{equation}
\left( \mathbf{P}_{0}\mathbf{P}_{1}.\mathbf{Q}_{0}\mathbf{Q}_{1}\right)
=\sigma \left( P_{0},Q_{1}\right) +\sigma \left( P_{1},Q_{0}\right) -\sigma
\left( P_{0},Q_{0}\right) -\sigma \left( P_{1},Q_{1}\right)  \label{b11a}
\end{equation}
Vectors $\mathbf{P}_{0}\mathbf{P}_{i}$, $\;i=1,2,...n$ are basic vectors of
the rectilinear coordinate system $K_{n}$ with the origin at the point $%
P_{0} $, and the metric tensors $g_{ik}\left( \mathcal{P}^{n}\right) $, $%
g^{ik}\left( \mathcal{P}^{n}\right) $, \ $i,k=1,2,...n$ in $K_{n}$ are
defined by the relations 
\begin{equation}
\sum\limits_{k=1}^{k=n}g^{ik}\left( \mathcal{P}^{n}\right) g_{lk}\left( 
\mathcal{P}^{n}\right) =\delta _{l}^{i},\qquad g_{il}\left( \mathcal{P}%
^{n}\right) =\left( \mathbf{P}_{0}\mathbf{P}_{i}.\mathbf{P}_{0}\mathbf{P}%
_{l}\right) ,\qquad i,l=1,2,...n  \label{a15b}
\end{equation}

III. Linear structure of the Euclidean space: 
\begin{equation}
\sigma \left( P,Q\right) =\frac{1}{2}\sum\limits_{i,k=1}^{i,k=n}g^{ik}\left( 
\mathcal{P}^{n}\right) \left( x_{i}\left( P\right) -x_{i}\left( Q\right)
\right) \left( x_{k}\left( P\right) -x_{k}\left( Q\right) \right) ,\qquad
\forall P,Q\in \Omega  \label{a15a}
\end{equation}
where coordinates $x_{i}\left( P\right) ,$ $i=1,2,...n$ of the point $P$ are
covariant coordinates of the vector $\mathbf{P}_{0}\mathbf{P}$, defined by
the relation 
\begin{equation}
x_{i}\left( P\right) =\left( \mathbf{P}_{0}\mathbf{P}_{i}.\mathbf{P}_{0}%
\mathbf{P}\right) ,\qquad i=1,2,...n  \label{b12}
\end{equation}

IV: The metric tensor matrix $g_{lk}\left( \mathcal{P}^{n}\right) $ has only
positive eigenvalues 
\begin{equation}
g_{k}>0,\qquad k=1,2,...,n  \label{a15c}
\end{equation}

V. The continuity condition: the system of equations 
\begin{equation}
\left( \mathbf{P}_{0}\mathbf{P}_{i}.\mathbf{P}_{0}\mathbf{P}\right)
=y_{i}\in \Bbb{R},\qquad i=1,2,...n  \label{b14}
\end{equation}
considered to be equations for determination of the point $P$ as a function
of coordinates $y=\left\{ y_{i}\right\} $,\ \ $i=1,2,...n$ has always one
and only one solution. Conditions II -- V contain a reference to the
dimension $n$ of the Euclidean space.

As far as the $\sigma $-immanent description of the proper Euclidean
geometry is possible, it is possible for any T-geometry, because any
geometrical object $\mathcal{O}$ and any relation $\mathcal{R}$ in the
physical geometry $\mathcal{G}$ is obtained from the corresponding
geometrical object $\mathcal{O}_{\mathrm{E}}$ and from the corresponding
relation $\mathcal{R}_{\mathrm{E}}$ in the proper Euclidean geometry $%
\mathcal{G}_{\mathrm{E}}$ by means of the replacement $\sigma _{\mathrm{E}%
}\rightarrow \sigma $ in description of $\mathcal{O}_{\mathrm{E}}$ an $%
\mathcal{R}_{\mathrm{E}}$. For such a replacement be possible, the
description of $\mathcal{O}_{\mathrm{E}}$ and $\mathcal{R}_{\mathrm{E}}$ is
not to refer to special properties of $\sigma _{\mathrm{E}}$, described by
conditions II -- V. A formal indicator of the conditions II -- V application
is a reference to the dimension $n$, because any of conditions II -- V
contains a reference to the dimension $n$ of the proper Euclidean space.

If nevertheless we use one of special properties II -- V of the Euclidean
space in the $\sigma $-immanent description of a geometrical object $%
\mathcal{O}$, or relation $\mathcal{R}$ , we refer to the dimension $n$ and,
ultimately, to the coordinate system, which is only a means of description.

Let us show this in the example of the determination of the straight in the $%
n$-dimensional Euclidean space. The straight $\mathcal{T}_{P_{0}Q}$ in the
proper Euclidean space is defined by two its points $P_{0}$ and $Q$ $%
\;\left( P_{0}\neq Q\right) $ as the set of points $R$ 
\begin{equation}
\mathcal{T}_{P_{0}Q}=\left\{ R\;|\;\mathbf{P}_{0}\mathbf{Q}||\mathbf{P}_{0}%
\mathbf{R}\right\}  \label{b15}
\end{equation}
where condition $\mathbf{P}_{0}\mathbf{Q}||\mathbf{P}_{0}\mathbf{R}$ means
that vectors $\mathbf{P}_{0}\mathbf{Q}$ and $\mathbf{P}_{0}\mathbf{R}$ are
collinear, i.e. the scalar product $\left( \mathbf{P}_{0}\mathbf{Q}.\mathbf{P%
}_{0}\mathbf{R}\right) $ of these two vectors satisfies the relation 
\begin{equation}
\mathbf{P}_{0}\mathbf{Q}||\mathbf{P}_{0}\mathbf{R:\qquad }\left( \mathbf{P}%
_{0}\mathbf{Q}.\mathbf{P}_{0}\mathbf{R}\right) ^{2}=\left( \mathbf{P}_{0}%
\mathbf{Q}.\mathbf{P}_{0}\mathbf{Q}\right) \left( \mathbf{P}_{0}\mathbf{R}.%
\mathbf{P}_{0}\mathbf{R}\right)  \label{b16}
\end{equation}
where the scalar product is defined by the relation (\ref{b11a}). Thus, the
straight line $\mathcal{T}_{P_{0}Q}$ is defined $\sigma $-immanently, i.e.
in terms of the world function $\sigma $. We shall use two different names
(straight and tube) for the geometric object $\mathcal{T}_{P_{0}Q}$. We
shall use the term ''straight'', when we want to stress that $\mathcal{T}%
_{P_{0}Q}$ is a result of deformation of the Euclidean straight. We shall
use the term ''tube'', when we want to stress that $\mathcal{T}_{P_{0}Q}$
may be a many-dimensional surface.

In the Euclidean geometry one can use another definition of collinearity.
Vectors $\mathbf{P}_{0}\mathbf{Q}$ and $\mathbf{P}_{0}\mathbf{R}$ are
collinear, if components of vectors $\mathbf{P}_{0}\mathbf{Q}$ and $\mathbf{P%
}_{0}\mathbf{R}$ in some coordinate system are proportional. For instance,
in the $n$-dimensional Euclidean space one can introduce rectangular
coordinate system, choosing $n+1$ points $\mathcal{P}^{n}=\left\{
P_{0},P_{1},...P_{n}\right\} $ and forming $n$ basic vectors $\mathbf{P}_{0}%
\mathbf{P}_{i}$, $i=1,2,...n$. Then the collinearity condition can be
written in the form of $n$ equations 
\begin{equation}
\mathbf{P}_{0}\mathbf{Q}||\mathbf{P}_{0}\mathbf{R:\qquad }\left( \mathbf{P}%
_{0}\mathbf{P}_{i}.\mathbf{P}_{0}\mathbf{Q}\right) =a\left( \mathbf{P}_{0}%
\mathbf{P}_{i}.\mathbf{P}_{0}\mathbf{R}\right) ,\qquad i=1,2,...n,\qquad
a\in \Bbb{R}  \label{b17}
\end{equation}
where $a$ is some constant. Relations (\ref{b17}) are relations for
covariant components of vectors $\mathbf{P}_{0}\mathbf{Q}$ and $\mathbf{P}%
_{0}\mathbf{R}$ in the considered coordinate system with basic vectors $%
\mathbf{P}_{0}\mathbf{P}_{i}$, $i=1,2,...n$. Let points $\mathcal{P}^{n}$ be
chosen in such a way, that $\left( \mathbf{P}_{0}\mathbf{P}_{1}.\mathbf{P}%
_{0}\mathbf{Q}\right) \neq 0$. Then eliminating the parameter $a$ from
relations (\ref{b17}), we obtain $n-1$ independent relations, and the
geometrical object 
\begin{eqnarray}
\mathcal{T}_{Q\mathcal{P}^{n}} &=&\left\{ R\;|\;\mathbf{P}_{0}\mathbf{Q}||%
\mathbf{P}_{0}\mathbf{R}\right\} =\bigcap\limits_{i=2}^{i=n}\mathcal{S}_{i},
\label{c2.1} \\
\mathcal{S}_{i} &=&\left\{ R\left| \frac{\left( \mathbf{P}_{0}\mathbf{P}_{i}.%
\mathbf{P}_{0}\mathbf{Q}\right) }{\left( \mathbf{P}_{0}\mathbf{P}_{1}.%
\mathbf{P}_{0}\mathbf{Q}\right) }=\frac{\left( \mathbf{P}_{0}\mathbf{P}_{i}.%
\mathbf{P}_{0}\mathbf{R}\right) }{\left( \mathbf{P}_{0}\mathbf{P}_{1}.%
\mathbf{P}_{0}\mathbf{R}\right) }\right. \right\} ,\qquad i=2,3,...n
\label{c2.2}
\end{eqnarray}
defined according to (\ref{b17}), depends on $n+2$ points $Q,\mathcal{P}^{n}$%
. This geometrical object $\mathcal{T}_{Q\mathcal{P}^{n}}$ is defined $%
\sigma $-immanently. It is a complex, consisting of the straight line and
the coordinate system, represented by $n+1$ points $\mathcal{P}^{n}=\left\{
P_{0},P_{1},...P_{n}\right\} $. In the Euclidean space the dependence on the
choice of the coordinate system and on $n+1$ points $\mathcal{P}^{n}$
determining this system, is fictitious. The geometrical object $\mathcal{T}%
_{Q\mathcal{P}^{n}}$ depends only on two points $P_{0},Q$ and coincides with
the straight line $\mathcal{T}_{P_{0}Q}$. But at deformations of the
Euclidean space the geometrical objects $\mathcal{T}_{Q\mathcal{P}^{n}}$ and 
$\mathcal{T}_{P_{0}Q}$ are deformed differently. The points $%
P_{1},P_{2},...P_{n}$ cease to be fictitious in definition of $\mathcal{T}_{Q%
\mathcal{P}^{n}}$, and geometrical objects $\mathcal{T}_{Q\mathcal{P}^{n}}$
and $\mathcal{T}_{P_{0}Q}$ become to be different geometric objects, in
general. But being different, in general, they may coincide in some special
cases.

What of the two geometrical objects in the deformed geometry should be
interpreted as the straight line, passing through the points $P_{0}$ and $Q$
in the geometry $\mathcal{G}$? Of course, it is $\mathcal{T}_{P_{0}Q}$,
because its definition does not contain a reference to a coordinate system,
whereas definition of $\mathcal{T}_{Q\mathcal{P}^{n}}$ depends on the choice
of the coordinate system, represented by points $\mathcal{P}^{n}$. In
general, definitions of geometric objects and relations between them are not
to refer to the means of description.

But in the given case the geometrical object $\mathcal{T}_{P_{0}Q}$ is, in
general, $(n-1)$-dimensional surface, whereas $\mathcal{T}_{Q\mathcal{P}%
^{n}} $ is an intersection of $(n-1)\;\;$ $(n-1)$-dimensional surfaces, i.e. 
$\mathcal{T}_{Q\mathcal{P}^{n}}$ is, in general, a one-dimensional curve.
The one-dimensional curve $\mathcal{T}_{Q\mathcal{P}^{n}}$ corresponds
better to our ideas on the straight line, than the $(n-1)$-dimensional
surface $\mathcal{T}_{P_{0}Q}$. Nevertheless, in physical geometry $\mathcal{%
G}$ it is $\mathcal{T}_{P_{0}Q}$, that is an analog of the Euclidean
straight line.

It is very difficult to overcome our conventional idea that the Euclidean
straight line cannot be deformed into many-dimensional surface, and \textit{%
this idea has been prevent for years from construction of T-geometries}.
Practically one uses such physical geometries, where deformation of the
Euclidean space transforms the Euclidean straight lines into one-dimensional
lines. It means that one chooses such geometries, where geometrical objects $%
\mathcal{T}_{P_{0}Q}$ and $\mathcal{T}_{Q\mathcal{P}^{n}}$ coincide. 
\begin{equation}
\mathcal{T}_{P_{0}Q}=\mathcal{T}_{Q\mathcal{P}^{n}}  \label{b19}
\end{equation}
Condition (\ref{b19}) of coincidence of the objects $\mathcal{T}_{P_{0}Q}$
and $\mathcal{T}_{Q\mathcal{P}^{n}}$, imposed on the T-geometry, restricts
list of possible T-geometries.

Let us consider the metric geometry, given on the set $\Omega $ of points.
The metric space $M=\left\{ \rho ,\Omega \right\} $ is given by the metric
(distance) $\rho $. 
\begin{eqnarray}
\rho &:&\quad \Omega \times \Omega \rightarrow \lbrack 0,\infty )\subset 
\Bbb{R}  \label{c2.3} \\
\rho (P,P) &=&0,\qquad \rho (P,Q)=\rho (Q,P),\qquad \forall P,Q\in \Omega
\label{c2.4} \\
\rho (P,Q) &\geq &0,\qquad \rho (P,Q)=0,\quad \text{iff }P=Q,\qquad \forall
P,Q\in \Omega  \label{c2.5} \\
0 &\leq &\rho (P,R)+\rho (R,Q)-\rho (P,Q),\qquad \forall P,Q,R\in \Omega
\label{c2.6}
\end{eqnarray}
At first sight the metric space is a special case of the $\sigma $-space (%
\ref{a1.3}), and the metric geometry is a special case of the T-geometry
with additional constraints (\ref{c2.5}), (\ref{c2.6}) imposed on the world
function $\sigma =\frac{1}{2}\rho ^{2}$. However it is not so, because the
metric geometry does not use the deformation principle. The fact, that the
Euclidean geometry can be described $\sigma $-immanently, as well as the
conditions (\ref{b10}) - (\ref{b14}), were not known until 1990. Additional
(with respect to the $\sigma $-space) constraints (\ref{c2.5}), (\ref{c2.6})
are imposed to eliminate the situation, when the straight line is not a
one-dimensional line. The fact is that, in the metric geometry the shortest
(straight) line can be constructed only in the case, when it is
one-dimensional.

Let us consider the set $\mathcal{EL}\left( P,Q,a\right) $ of points $R$ 
\begin{equation}
\mathcal{EL}\left( P,Q,a\right) =\left\{ R|f_{P,Q,a}\left( R\right)
=0\right\} ,\qquad f_{P,Q,a}\left( R\right) =\rho (P,R)+\rho (R,Q)-2a
\label{c2.8}
\end{equation}
If the metric space coincides with the proper Euclidean space, this set of
points is an ellipsoid with focuses at the points $P,Q$ and the large
semiaxis $a$. The relations $f_{P,Q,a}\left( R\right) >0$, $f_{P,Q,a}\left(
R\right) =0$, $f_{P,Q,a}\left( R\right) <0$ determine respectively external
points, boundary points and internal points of the ellipsoid. If $\rho
\left( P,Q\right) =2a$, we obtain the degenerate ellipsoid, which coincides
with the segment $\mathcal{T}_{\left[ PQ\right] }$ of the straight line,
passing through the points $P$, $Q$.\ In the proper Euclidean geometry, the
degenerate ellipsoid is one-dimensional segment of the straight line, but it
is not evident that it is one-dimensional in the case of arbitrary metric
geometry. For such a degenerate ellipsoid be one-dimensional in the
arbitrary metric space, it is necessary that any degenerate ellipsoid $%
\mathcal{EL}\left( P,Q,\rho \left( P,Q\right) /2\right) $ have no internal
points. This constraint is written in the form 
\begin{equation}
f_{P,Q,\rho \left( P,Q\right) /2}\left( R\right) =\rho (P,R)+\rho (R,Q)-\rho
(P,Q)\geq 0  \label{c2.9}
\end{equation}

Comparing relation (\ref{c2.9}) with (\ref{c2.6}), we see that the
constraint (\ref{c2.6}) is introduced to make the straight (shortest) line
to be one-dimensional (absence of internal points in the geometrical object
determined by two points).

As far as the metric geometry does not use the deformation principle, it is
a poor geometry, because in the framework of this geometry one cannot
construct the scalar product of two vectors, define linear independence of
vectors and construct such geometrical objects as planes. All these objects
as well as other are constructed on the basis of the deformation of the
proper Euclidean geometry.

Generalizing the metric geometry, Menger \cite{M28} and Blumenthal \cite{B53}
removed the triangle axiom (\ref{c2.6}). They tried to construct the
distance geometry, which would be a more general geometry, than the metric
one. As far as they did not use the deformation principle, they could not
determine the shortest (straight) line without a reference to the
topological concept of the curve $\mathcal{L}$, defined as a continuous
mapping 
\begin{equation}
\mathcal{L}:\qquad \left[ 0,1\right] \rightarrow \Omega  \label{a1.1}
\end{equation}
which cannot be expressed only via the distance. As a result the distance
geometry appeared to be not a pure metric geometry, what the T-geometry is.

\section{Conditions of the deformation principle \newline
application}

Riemannian geometries satisfy the condition (\ref{b19}). The Riemannian
geometry is a kind of inhomogeneous physical geometry, and, hence, it uses
the deformation principle. Constructing the Riemannian geometry, the
infinitesimal Euclidean distance is deformed into the Riemannian distance.
The deformation is chosen in such a way that any Euclidean straight line $%
\mathcal{T}_{\mathrm{E}P_{0}Q}$, passing through the point $P_{0}$,
collinear to the vector $\mathbf{P}_{0}\mathbf{Q}$, transforms into the
geodesic $\mathcal{T}_{P_{0}Q}$, passing through the point $P_{0}$,
collinear to the vector $\mathbf{P}_{0}\mathbf{Q}$ in the Riemannian space.

Note that in T-geometries, satisfying the condition (\ref{b19}) for all
points $Q,\mathcal{P}^{n}$, the straight line 
\begin{equation}
\mathcal{T}_{Q_{0};P_{0}Q}=\left\{ R\;|\;\mathbf{P}_{0}\mathbf{Q}||\mathbf{Q}%
_{0}\mathbf{R}\right\}  \label{b3.0}
\end{equation}
passing through the point $Q_{0}$ collinear to the vector $\mathbf{P}_{0}%
\mathbf{Q}$, is not a one-dimensional line, in general. If the Riemannian
geometries be T-geometries, they would contain non-one-dimensional geodesics
(straight lines). But the Riemannian geometries are not T-geometries,
because at their construction one uses not only the deformation principle,
but some other methods, containing a reference to the means of description.
In particular, in the Riemannian geometries the absolute parallelism is
absent, and one cannot to define a straight line (\ref{b3.0}), because the
relation $\mathbf{P}_{0}\mathbf{Q}||\mathbf{Q}_{0}\mathbf{R}$ is not
defined, if points $P_{0}$ and $Q_{0}$ do not coincide. On one hand, a lack
of absolute parallelism allows one to go around the problem of
non-one-dimensional straight lines. On the other hand, it makes the
Riemannian geometries to be inconsistent, because they cease to be
T-geometries, which are consistent by the construction (see for details \cite
{R02}).

The fact is that the application of \textit{only deformation principle }is
sufficient for construction of a physical geometry. Besides, such a
construction is consistent, because the original Euclidean geometry is
consistent and, deforming it, we do not use any reasonings. If we introduce
additional structure (for instance, a topological structure) we obtain a
fortified physical geometry, i.e. a physical geometry with additional
structure on it. The physical geometry with additional structure on it is a
more pithy construction, than the physical geometry simply. But it is valid
only in the case, when we consider the additional structure as an addition
to the physical geometry. If we use an additional structure in construction
of the geometry, we identify the additional structure with one of structures
of the physical geometry. If we demand that the additional structure to be a
structure of physical geometry, we restrict an application of the
deformation principle and reduce the list of possible physical geometries,
because coincidence of the additional structure with some structure of a
physical geometry is possible not for all physical geometries, but only for
some of them.

Let, for instance, we use concept of a curve $\mathcal{L}$ (\ref{a1.1}) for
construction of a physical geometry. The concept of curve $\mathcal{L}$,
considered as a continuous mapping is a topological structure, which cannot
be expressed only via the distance or via the world function. A use of the
mapping (\ref{a1.1}) needs an introduction of topological space and, in
particular, the concept of continuity. If we identify the topological curve (%
\ref{a1.1}) with the ''metrical'' curve, defined as a broken line 
\begin{equation}
\mathcal{T}_{\mathrm{br}}=\bigcup\limits_{i}\mathcal{T}_{\left[ P_{i}P_{i+1}%
\right] },\qquad \mathcal{T}_{\left[ P_{i}P_{i+1}\right] }=\left\{ R|\sqrt{%
2\sigma \left( P_{i},P_{i+1}\right) }-\sqrt{2\sigma \left( P_{i},R\right) }-%
\sqrt{2\sigma \left( R,P_{i+1}\right) }\right\}  \label{a1.2}
\end{equation}
consisting of the straight line segments $\mathcal{T}_{\left[ P_{i}P_{i+1}%
\right] }$ between the points $P_{i}$, $P_{i+1}$, we truncate the list of
possible geometries, because such an identification is possible only in some
physical geometries. Identifying (\ref{a1.1}) and (\ref{a1.2}), we eliminate
all discrete physical geometries and those continuous physical geometries,
where the segment $\mathcal{T}_{\left[ P_{i}P_{i+1}\right] }$ of straight
line is a surface, but not a one-dimensional set of points. Thus, additional
structures may lead to (i) a fortified physical geometry, (ii) a restricted
physical geometry and (iii) a restricted fortified physical geometry. The
result depends on the method of the additional structure application.

Note that some constraints (continuity, convexity, lack of absolute
parallelism), imposed on physical geometries are a result of a disagreement
of the applied means of the geometry construction. In the T-geometry, which
uses only the deformation principle, there is no such restrictions. Besides,
the T-geometry accepts some new property of a physical geometry, which is
not accepted by conventional versions of physical geometry. This property,
called the geometry nondegeneracy, follows directly from the application of
arbitrary deformations to the proper Euclidean geometry.

The geometry is degenerate at the point $P_{0}$ in the direction of the
vector $\mathbf{Q}_{0}\mathbf{Q}$, $\left| \mathbf{Q}_{0}\mathbf{Q}\right|
\neq 0$, if the relations 
\begin{equation}
\mathbf{Q}_{0}\mathbf{Q}\uparrow \uparrow \mathbf{P}_{0}\mathbf{R:\qquad }%
\left( \mathbf{Q}_{0}\mathbf{Q}.\mathbf{P}_{0}\mathbf{R}\right) =\sqrt{%
\left| \mathbf{Q}_{0}\mathbf{Q}\right| \cdot \left| \mathbf{P}_{0}\mathbf{R}%
\right| },\qquad \left| \mathbf{P}_{0}\mathbf{R}\right| =a\neq 0
\label{b3.1}
\end{equation}
considered as equations for determination of the point $R$, have not more,
than one solution for any $a\neq 0$. Otherwise, the geometry is
nondegenerate at the point $P_{0}$ in the direction of the vector $\mathbf{Q}%
_{0}\mathbf{Q}$. Note that the first equation (\ref{b3.1}) is the condition
of the parallelism of vectors $\mathbf{Q}_{0}\mathbf{Q}$ and $\mathbf{P}_{0}%
\mathbf{R}$.

The proper Euclidean geometry is degenerate, i.e. it is degenerate at all
points in directions of all vectors. Considering the Minkowski geometry, one
should distinguish between the Minkowski T-geometry and Minkowski geometry.
The two geometries are described by the same world function and differ in
the definition of the parallelism. In the Minkowski T-geometry the
parallelism of two vectors $\mathbf{\mathbf{Q}_{0}\mathbf{Q}}$ and $\mathbf{%
\mathbf{P}_{0}\mathbf{R}}$ is defined by the first equation (\ref{b3.1}).
This definition is based on the deformation principle. In Minkowski geometry
the parallelism is defined by the relation of the type of (\ref{b17}) 
\begin{equation}
\mathbf{Q}_{0}\mathbf{Q}\uparrow \uparrow \mathbf{P}_{0}\mathbf{R:\qquad }%
\left( \mathbf{P}_{0}\mathbf{P}_{i}.\mathbf{Q}_{0}\mathbf{Q}\right) =a\left( 
\mathbf{P}_{0}\mathbf{P}_{i}.\mathbf{P}_{0}\mathbf{R}\right) ,\qquad
i=1,2,...n,\qquad a>0  \label{b3.1a}
\end{equation}
where points $\mathcal{P}^{n}=\left\{ P_{0},P_{1},...P_{n}\right\} $
determine a rectilinear coordinate system with basic vectors $\mathbf{P}_{0}%
\mathbf{P}_{i}$, $i=1,2,..n$ in the $n$-dimensional Minkowski geometry ($n$%
-dimensional pseudo-Euclidean geometry of index $1$). Dependence of the
definition (\ref{b3.1a}) on the points $\left( P_{1},P_{2},...P_{n}\right) $
is fictitious, but dependence on the number $n+1$ of points $\mathcal{P}^{n}$
is essential. Thus, definition (\ref{b3.1a}) depends on the method of the
geometry description.

The Minkowski T-geometry is degenerate at all points in direction of all
timelike vectors, and it is nondegenerate at all points in direction of all
spacelike vectors. The Minkowski geometry is degenerate at all points in
direction of all vectors. Conventionally one uses the Minkowski geometry,
ignoring the nondegeneracy in spacelike directions.

Considering the proper Riemannian geometry, one should distinguish between
the Riemannian T-geometry and the Riemannian geometry. The two geometries
are described by the same world function. They differ in the definition of
the parallelism. In the Riemannian T-geometry the parallelism of two vectors 
$\mathbf{\mathbf{Q}_{0}\mathbf{Q}}$ and $\mathbf{\mathbf{P}_{0}\mathbf{R}}$
is defined by the first equation (\ref{b3.1}). In the Riemannian geometry
the parallelism of two vectors $\mathbf{\mathbf{Q}_{0}\mathbf{Q}}$ and $%
\mathbf{\mathbf{P}_{0}\mathbf{R}}$ is defined only in the case, when the
points $P_{0}$ and $Q_{0}$ coincide. Parallelism of remote vectors $\mathbf{%
\mathbf{Q}_{0}\mathbf{Q}}$ and $\mathbf{\mathbf{P}_{0}\mathbf{R}}$ is not
defined, in general. This fact is known as absence of absolute parallelism.

The proper Riemannian T-geometry is locally degenerate, i.e. it is
degenerate at all points $P_{0}$ in direction of vectors $\mathbf{P}_{0}%
\mathbf{Q}$. In the general case, when $P_{0}\neq Q_{0}$, the proper
Riemannian T-geometry is nondegenerate, in general. The proper Riemannian
geometry is degenerate, because it is degenerate locally, whereas the
nonlocal degeneracy is not defined in the Riemannian geometry, because of
the lack of absolute parallelism. Conventionally one uses the Riemannian
geometry (not Rienannian T-geometry) and ignores the property of the
nondegenracy completely.

From the viewpoint of the conventional approach to the physical geometry the
nondegeneracy is an undesirable property of a physical geometry, although
from the logical viewpoint and from viewpoint of the deformation principle
the nondegenracy is an inherent property of a physical geometry. The
nonlocal nondegeneracy is ejected from the proper Riemannian geometry by
denial of existence of the remote vector parallelism. Nondegeneracy in the
spacelike directions is ejected from the Minkowski geometry by means of the
redefinition of the two vectors parallelism. But the nondegeneracy is an
important property of the real space-time geometry. To appreciate this, let
us consider an example.

\section{Simple example of nondegenerate space-time geometry}

The T-geometry \cite{R01} is defined on the $\sigma $-space $V=\left\{
\sigma ,\Omega \right\} $, where $\Omega $ is an arbitrary set of points and
the world function $\sigma $ is defined by the relations

\begin{equation}
\sigma :\qquad \Omega \times \Omega \rightarrow \Bbb{R},\qquad \sigma \left(
P,Q\right) =\sigma \left( Q,P\right) ,\qquad \sigma \left( P,P\right)
=0,\qquad \forall P,Q\in \Omega  \label{b3.2}
\end{equation}
Geometrical objects (vector $\mathbf{PQ}$, scalar product of vectors $\left( 
\mathbf{P}_{0}\mathbf{P}_{1}.\mathbf{Q}_{0}\mathbf{Q}_{1}\right) $,
collinearity of vectors $\mathbf{\ P}_{0}\mathbf{P}_{1}||\mathbf{Q}_{0}%
\mathbf{Q}_{1}$, segment of straight line $\mathcal{T}_{\left[ P_{0}P_{1}%
\right] }$, etc.) are defined on the $\sigma $-space in the same way, as
they are defined $\sigma $-immanently in the proper Euclidean space.
Practically one uses the deformation principle, although it is not mentioned
in all definitions.

Let us consider a simple example of the space-time geometry $\mathcal{G}_{%
\mathrm{d}}$, described by the T-geometry on 4-dimensional manifold $%
\mathcal{M}_{1+3}$. The world function $\sigma _{\mathrm{d}}$ is described
by the relation 
\begin{equation}
\sigma _{\mathrm{d}}=\sigma _{\mathrm{M}}+D\left( \sigma _{\mathrm{M}%
}\right) ,\qquad D\left( \sigma _{\mathrm{M}}\right) =\left\{ 
\begin{array}{ll}
\sigma _{\mathrm{M}}+d & \text{if\ }\sigma _{0}<\sigma _{\mathrm{M}} \\ 
\left( 1+\frac{d}{\sigma _{0}}\right) \sigma _{\mathrm{M}} & \text{if\ }%
0\leq \sigma _{\mathrm{M}}\leq \sigma _{0} \\ 
\sigma _{\mathrm{M}} & \text{if\ }\sigma _{\mathrm{M}}<0
\end{array}
\right.  \label{b3.3}
\end{equation}
where $d\geq 0$ and $\sigma _{0}>0$ are some constants. The quantity $\sigma
_{\mathrm{M}}$ is the world function in the Minkowski space-time geometry $%
\mathcal{G}_{\mathrm{M}}$. In the orthogonal rectilinear (inertial)
coordinate system $x=\left( t,\mathbf{x}\right) $ the world function $\sigma
_{\mathrm{M}}$ has the form 
\begin{equation}
\sigma _{\mathrm{M}}\left( x,x^{\prime }\right) =\frac{1}{2}\left(
c^{2}\left( t-t^{\prime }\right) ^{2}-\left( \mathbf{x}-\mathbf{x}^{\prime
}\right) ^{2}\right)  \label{b3.4}
\end{equation}
where $c$ is the speed of the light.

Let us compare the broken line (\ref{a1.2}) in Minkowski space-time geometry 
$\mathcal{G}_{\mathrm{M}}$ and in the distorted geometry $\mathcal{G}_{%
\mathrm{d}}$. We suppose that $\mathcal{T}_{\mathrm{br}}$ is timelike broken
line, and all links $\mathcal{T}_{\left[ P_{i}P_{i+1}\right] }$ of $\mathcal{%
T}_{\mathrm{br}}$ are timelike and have the same length 
\begin{equation}
\left| \mathbf{P}_{i}\mathbf{P}_{i+1}\right| _{\mathrm{d}}=\sqrt{2\sigma _{%
\mathrm{d}}\left( P_{i},P_{i+1}\right) }=\mu _{\mathrm{d}}>0,\qquad i=0,\pm
1,\pm 2,...  \label{b3.5}
\end{equation}
\begin{equation}
\left| \mathbf{P}_{i}\mathbf{P}_{i+1}\right| _{\mathrm{M}}=\sqrt{2\sigma _{%
\mathrm{M}}\left( P_{i},P_{i+1}\right) }=\mu _{\mathrm{M}}>0,\qquad i=0,\pm
1,\pm 2,...  \label{b3.5a}
\end{equation}
where indices ''d'' and ''M'' mean that the quantity is calculated by means
of $\sigma _{\mathrm{d}}$ and $\sigma _{\mathrm{M}}$ respectively. Vector $%
\mathbf{P}_{i}\mathbf{P}_{i+1}$ is regarded as the momentum of the particle
at the segment $\mathcal{T}_{\left[ P_{i}P_{i+1}\right] }$, and the quantity 
$\left| \mathbf{P}_{i}\mathbf{P}_{i+1}\right| =\mu $ is interpreted as its
(geometric) mass. It follows from definition (\ref{b11a}) and relation (\ref
{b3.3}), that for timelike vectors $\mathbf{P}_{i}\mathbf{P}_{i+1}$ with $%
\mu >\sqrt{2\sigma _{0}}$%
\begin{equation}
\left| \mathbf{P}_{i}\mathbf{P}_{i+1}\right| _{\mathrm{d}}^{2}=\mu _{\mathrm{%
d}}^{2}=\mu _{\mathrm{M}}^{2}+2d,\qquad \mu _{\mathrm{M}}^{2}>2\sigma _{0}
\label{b3.6}
\end{equation}
\begin{equation}
\left( \mathbf{P}_{i-1}\mathbf{P}_{i}.\mathbf{P}_{i}\mathbf{P}_{i+1}\right)
_{\mathrm{d}}=\left( \mathbf{P}_{i-1}\mathbf{P}_{i}.\mathbf{P}_{i}\mathbf{P}%
_{i+1}\right) _{\mathrm{M}}+d  \label{b3.7}
\end{equation}
Calculation of the shape of the segment $\mathcal{T}_{\left[ P_{0}P_{1}%
\right] }\left( \sigma _{\mathrm{d}}\right) $ in $\mathcal{G}_{\mathrm{d}}$
gives the relation 
\begin{equation}
r^{2}(\tau )=\left\{ 
\begin{array}{ll}
\tau ^{2}\mu _{\mathrm{d}}^{2}\frac{\left( 1-\frac{\tau d}{2\left( \sigma
_{0}+d\right) }\right) ^{2}}{\left( 1-\frac{2d}{\mu _{\mathrm{d}}^{2}}%
\right) }-\frac{\tau ^{2}\mu _{\mathrm{d}}^{2}\sigma _{0}}{\left( \sigma
_{0}+d\right) }, & 0<\tau <\frac{\sqrt{2(\sigma _{0}+d)}}{\mu _{\mathrm{d}}}
\\ 
\frac{3d}{2}+2d\left( \tau -1/2\right) ^{2}\left( 1-\frac{2d}{\mu _{\mathrm{d%
}}^{2}}\right) ^{-1}, & \frac{\sqrt{2(\sigma _{0}+d)}}{\mu _{\mathrm{d}}}%
<\tau <1-\frac{\sqrt{2(\sigma _{0}+d)}}{\mu _{\mathrm{d}}} \\ 
\left( 1-\tau \right) ^{2}\mu _{\mathrm{d}}^{2}\left[ \frac{\left( 1-\frac{%
\left( 1-\tau \right) d}{2\left( \sigma _{0}+d\right) }\right) ^{2}}{\left(
1-\frac{2d}{\mu _{\mathrm{d}}^{2}}\right) }-\frac{\sigma _{0}}{\left( \sigma
_{0}+d\right) }\right] , & 1-\frac{\sqrt{2(\sigma _{0}+d)}}{\mu _{\mathrm{d}}%
}<\tau <1
\end{array}
\right. ,  \label{b3.7a}
\end{equation}
where $r\left( \tau \right) $ is the spatial radius of the segment $\mathcal{%
T}_{\left[ P_{0}P_{1}\right] }\left( \sigma _{\mathrm{d}}\right) $ in the
coordinate system, where points $P_{0}$ and $P_{1}$ have coordinates $%
P_{0}=\left\{ 0,0,0,0\right\} $, $P_{1}=\left\{ \mu _{\mathrm{d}%
},0,0,0\right\} $ and $\tau $ is a parameter along the segment $\mathcal{T}_{%
\left[ P_{0}P_{1}\right] }\left( \sigma _{\mathrm{d}}\right) $ ($\tau \left(
P_{0}\right) =0$, $\tau \left( P_{1}\right) =1$). One can see from (\ref
{b3.7a}) that the characteristic value of the segment radius is $\sqrt{d}$.

Let the broken tube $\mathcal{T}_{\mathrm{br}}$ describe the ''world line''
of a free particle. It means by definition that any link $\mathbf{P}_{i-1}%
\mathbf{P}_{i}$ is parallel to the adjacent link $\mathbf{P}_{i}\mathbf{P}%
_{i+1}$%
\begin{equation}
\mathbf{P}_{i-1}\mathbf{P}_{i}\uparrow \uparrow \mathbf{P}_{i}\mathbf{P}%
_{i+1}:\qquad \left( \mathbf{P}_{i-1}\mathbf{P}_{i}.\mathbf{P}_{i}\mathbf{P}%
_{i+1}\right) -\left| \mathbf{P}_{i-1}\mathbf{P}_{i}\right| \cdot \left| 
\mathbf{P}_{i}\mathbf{P}_{i+1}\right| =0  \label{b3.8}
\end{equation}
Definition of parallelism is different in geometries $\mathcal{G}_{\mathrm{M}%
}$ and $\mathcal{G}_{\mathrm{d}}$. As a result links, which are parallel in
the geometry $\mathcal{G}_{\mathrm{M}}$, are not parallel in $\mathcal{G}_{%
\mathrm{d}}$ and vice versa.

Let $\mathcal{T}_{\mathrm{br}}\left( \sigma _{\mathrm{M}}\right) $ describe
the world line of a free particle in the geometry $\mathcal{G}_{\mathrm{M}}$%
. The angle $\vartheta _{\mathrm{M}}$ between the adjacent links in $%
\mathcal{G}_{\mathrm{M}}$ is defined by the relation 
\begin{equation}
\cosh \vartheta _{\mathrm{M}}=\frac{\left( \mathbf{P}_{-1}\mathbf{P}_{0}.%
\mathbf{P}_{0}\mathbf{P}_{1}\right) _{\mathrm{M}}}{\left| \mathbf{P}_{0}%
\mathbf{P}_{1}\right| _{\mathrm{M}}\cdot \left| \mathbf{P}_{-1}\mathbf{P}%
_{0}\right| _{\mathrm{M}}}=1  \label{b3.9}
\end{equation}
The angle $\vartheta _{\mathrm{M}}=0$, and the geometrical object $\mathcal{T%
}_{\mathrm{br}}\left( \sigma _{\mathrm{M}}\right) $ is a timelike straight
line on the manifold $\mathcal{M}_{1+3}$.

Let now $\mathcal{T}_{\mathrm{br}}\left( \sigma _{\mathrm{d}}\right) $
describe the world line of a free particle in the geometry $\mathcal{G}_{%
\mathrm{d}}$. The angle $\vartheta _{\mathrm{d}}$ between the adjacent links
in $\mathcal{G}_{\mathrm{d}}$ is defined by the relation 
\begin{equation}
\cosh \vartheta _{\mathrm{d}}=\frac{\left( \mathbf{P}_{i-1}\mathbf{P}_{i}.%
\mathbf{P}_{i}\mathbf{P}_{i+1}\right) _{\mathrm{d}}}{\left| \mathbf{P}_{i}%
\mathbf{P}_{i+1}\right| _{\mathrm{d}}\cdot \left| \mathbf{P}_{i-1}\mathbf{P}%
_{i}\right| _{\mathrm{d}}}=1  \label{b3.10}
\end{equation}
The angle $\vartheta _{\mathrm{d}}=0$ also. If we draw the broken tube $%
\mathcal{T}_{\mathrm{br}}\left( \sigma _{\mathrm{d}}\right) $ on the
manifold $\mathcal{M}_{1+3}$, using coordinates of basic points $P_{i}$ and
measure the angle $\vartheta _{\mathrm{dM}}$ between the adjacent links in
the Minkowski geometry $\mathcal{G}_{\mathrm{M}}$, we obtain for the angle $%
\vartheta _{\mathrm{dM}}$ the following relation 
\begin{equation}
\cosh \vartheta _{\mathrm{dM}}=\frac{\left( \mathbf{P}_{i-1}\mathbf{P}_{i}.%
\mathbf{P}_{i}\mathbf{P}_{i+1}\right) _{\mathrm{M}}}{\left| \mathbf{P}_{i}%
\mathbf{P}_{i+1}\right| _{\mathrm{M}}\cdot \left| \mathbf{P}_{i-1}\mathbf{P}%
_{i}\right| _{\mathrm{M}}}=\frac{\left( \mathbf{P}_{i-1}\mathbf{P}_{i}.%
\mathbf{P}_{i}\mathbf{P}_{i+1}\right) _{\mathrm{d}}-d}{\left| \mathbf{P}_{i}%
\mathbf{P}_{i+1}\right| _{\mathrm{d}}^{2}-2d}  \label{b3.11}
\end{equation}
Substituting the value of $\left( \mathbf{P}_{i-1}\mathbf{P}_{i}.\mathbf{P}%
_{i}\mathbf{P}_{i+1}\right) _{\mathrm{d}}$, taken from (\ref{b3.10}), we
obtain 
\begin{equation}
\cosh \vartheta _{\mathrm{dM}}=\frac{\mu _{\mathrm{d}}^{d}-d}{\mu _{\mathrm{d%
}}^{2}-2d}\approx 1+\frac{d}{\mu _{\mathrm{d}}^{2}},\qquad d\ll \mu _{%
\mathrm{d}}^{2}  \label{b3.12}
\end{equation}
Hence, $\vartheta _{\mathrm{dM}}\approx \sqrt{2d}/\mu _{\mathrm{d}}$. It
means, that the adjacent link is located on the cone of angle $\sqrt{2d}/\mu
_{\mathrm{d}}$, and the whole line $\mathcal{T}_{\mathrm{br}}\left( \sigma _{%
\mathrm{d}}\right) $ has a random shape, because any link wabbles with the
characteristic angle $\sqrt{2d}/\mu _{\mathrm{d}}$. The wabble angle depends
on the space-time distortion $d$ and on the particle mass $\mu _{\mathrm{d}}$%
. The wabble angle is small for the large mass of a particle. The random
displacement of the segment end is of the order $\mu _{\mathrm{d}}\vartheta
_{\mathrm{dM}}=\sqrt{2d}$, i.e. of the same order as the segment width. It
is reasonable, because these two phenomena have the common source: the
space-time distortion $D$.

One should note that the space-time geometry influences the stochasticity of
particle motion nonlocally in the sense, that the form of the world function
(\ref{b3.3}) for values of $\sigma _{\mathrm{M}}<\frac{1}{2}\mu _{\mathrm{d}%
}^{2}$ is unessential for the motion stochasticity of the particle of the
mass $\mu _{\mathrm{d}}$.

Such a situation, when the world line of a free particle is stochastic in
the deterministic geometry, and this stochasticity depends on the particle
mass, seems to be rather exotic and incredible. But experiments show that
the motion of real particles of small mass is stochastic indeed, and this
stochasticity increases, when the particle mass decreases. From physical
viewpoint a theoretical foundation of the stochasticity is desirable, and
some researchers invent stochastic geometries, noncommutative geometries and
other exotic geometrical constructions, to obtain the quantum stochasticity.
But in the Riemannian space-time geometry the particle motion does not
depend on the particle mass, and in the framework of the Riemannian
space-time geometry it is difficult to explain the quantum stochasticity by
the space-time geometry properties. Distorted geometry $\mathcal{G}_{\mathrm{%
d}}$ explains the stochasticity and its dependence on the particle mass
freely. Besides, at proper choice of the distortion $d$ the statistical
description of stochastic $\mathcal{T}_{\mathrm{br}}$ leads to the quantum
description (Schr\"{o}dinger equation) \cite{R91}. It is sufficient to set $%
d=0.5\hbar \left( bc\right) ^{-1}$, where $\hbar $ is the quantum constant, $%
c$ is the speed of the light, and $b$ is some universal constant, connecting
the geometrical mass $\mu $ with the usual particle mass $m$ by means of the
relation $m=b\mu $. In other words, the distorted space-time geometry (\ref
{b3.3}) is closer to the real space-time geometry, than the Minkowski
geometry $\mathcal{G}_{\mathrm{M}}$.

Further development of the statistical description of geometrical
stochasticity leads to a creation of the model conception of quantum
phenomena (MCQP), which relates to the conventional quantum theory
approximately in the same way as the statistical physics relates to the
axiomatic thermodynamics. MCQP is the well defined relativistic conception
with effective methods of investigation \cite{R03}, whereas the conventional
quantum theory is not well defined, because it uses incorrect space-time
geometry, whose incorrectness is compensated by additional hypotheses
(quantum principles). Besides, it has problems with application of the
nonrelativistic quantum mechanical technique to the description of
relativistic phenomena.

The geometry $\mathcal{G}_{\mathrm{d}}$ is a homogeneous geometry as well as
the Minkowski geometry, because the world function $\sigma _{\mathrm{d}}$ is
invariant with respect to all coordinate transformations, with respect to
which the world function $\sigma _{\mathrm{M}}$ is invariant. In this
connection the question arises, whether one could invent some axiomatics for 
$\mathcal{G}_{\mathrm{d}}$ and derive the geometry $\mathcal{G}_{\mathrm{d}}$
from this axiomatics by means of proper reasonings. Note that such an
axiomatics is to depend on the parameter $d$, because the world function $%
\sigma _{\mathrm{d}}$ depends on this parameter. If $d=0$, this axiomatics
is to coincide with the axiomatics of the Minkowski geometry $\mathcal{G}_{%
\mathrm{M}}$. If $d\neq 0$, this axiomatics cannot coincide with the
axiomatics of $\mathcal{G}_{\mathrm{M}}$, because some axioms of $\mathcal{G}%
_{\mathrm{M}}$ are not satisfied in this case. In general, the invention of
axiomatics, depending on the parameter $d$ and in the general case on the
distortion function $D$, seems to be a very difficult problem. Besides, why
invent the axiomatics? We had derived the axiomatics for the proper
Euclidean geometry, when we constructed it before. There is no necessity to
repeat this process any time, when we construct a new geometry. It is
sufficient to apply the deformation principle to the constructed Euclidean
geometry written $\sigma $-immanently. Application of the deformation
principle to the Euclidean geometry is a very simple and very general
procedue, which is not restricted by continuity, convexity and other
artificial constraints, generated by our preconceived approach to the
physical geometry. (Bias of the approach is displayed in the antecedent
supposition on the one-dimensionality of any straight line in any physical
geometry).

Thus, we have seen that the nondegeneracy of the physical geometry as well
as non-one-dimensionality of the straight line are properties of the real
physical geometries. The proper Euclidean geometry is a ground for all
physical geometries, and it is a degenerate geometry. Nevertheless, it is
beyond reason to deny an existence of nondegenerate physical geometries.

Thus, the deformation principle together with the $\sigma $-immanent
description appears to be a very effective mathematical tool for
construction of physical geometries.

\begin{enumerate}
\item  The deformation principle uses results obtained at construction of
the proper Euclidean geometry and does not add any additional supposition on
properties of geometrical objects.

\item  The deformation principle uses only the real characteristic of the
physical geometry -- its world function and does not use any additional
means of description.

\item  The deformation principle is very simple and allows one to
investigate only that part of geometry which one is interested in.

\item  Application of the deformation principle allows one to obtain the
true space-time geometry, whose unexpected properties cannot be obtained at
the conventional approach to physical geometry.
\end{enumerate}

\end{document}